\documentclass[11pt,a4paper]{article}
\usepackage{indentfirst}
\usepackage{amsmath,amssymb,amsfonts,amsthm,graphics}
\usepackage{makeidx}
\usepackage{color}
\usepackage[dvips]{graphicx}
\usepackage{eufrak}
\usepackage{mathrsfs}

\begin{document}
\title{\Large \bf Formally self-dual linear binary codes from circulant graphs\footnote{Supported by ``973" program No.2013CB834204.}}
\author{\small Xueliang Li, Yaping Mao, Meiqin Wei\\
\small Nankai University, Tianjin 300071, China\\
\small Center for Combinatorics and LPMC-TJKLC\\
\small E-mails: lxl@nankai.edu.cn; maoyaping@ymail.com; weimeiqin@mail.nankai.edu.cn\\
\small and\\
\small Ruihu Li\\
\small The air force engineering University\\
\small Institute of science, Xi'an 710051, China\\
\small E-mail: liruihu@aliyun.com}
\date{}
\maketitle
\begin{abstract}
In $2002$, Tonchev first constructed some linear binary codes
defined by the adjacency matrices of undirected graphs. So, graph is
an important tool for searching optimum codes. In this paper, we
introduce a new method of searching (proposed) optimum formally
self-dual linear binary codes from circulant graphs.\\[3mm]
{\bf AMS Subject Classification 2010}: 94B05, 05C50, 05C25.

\end{abstract}

\section{Introduction}

A \emph{linear binary code} $\mathcal {C}$ of length $n$ and
dimension $k$ (or an $[n,k]$ code), is a $k$-dimensional linear
subspace of the $n$-dimensional binary vector space $F^n_2$. The
\emph{Hamming distance} between two vectors $x=(x_1,\cdots,x_n)$,
$y=(y_1,\cdots,y_n)$ is equal to the number of indices $i$ such that
$x_i\neq y_i$. The \emph{Hamming weight} of a vector $x$, which is
denoted by $wt(x)$, is the number of its nonzero coordinates. The
\emph{minimum distance $d$} of a code is defined as the smallest
possible distance between pairs of distinct codewords. A generator
matrix for an $[n,k]$ code $\mathcal {C}$ is any $k\times n$ matrix
$G$ whose rows form a basis for $\mathcal {C}$. In general, there
are many generator matrices for a code.

We say that a code is optimum if it meets the lower and upper bounds
in the Code Tables, and a proposed optimum code if it only meets the
lower bound in the Code Tables. The distribution of a code is the
sequence $(A_{0},A_{1},\cdots,A_{n})$, where $A_{i}$ is the number
of codewords of weight $i$. The \emph{weight enumerator} of the code
is the polynomial
$$
W(z)=\sum\limits_{i=0}^{n}A_{i}z^{i}.
$$

Let us now introduce some concepts and notions from Graph Theory. An
\emph{undirected graph} $\Gamma=(V,E)$ is a set
$V(\Gamma)=\{v_1,v_2,\cdots,v_n\}$ of vertices together with a
collection $E(\Gamma)$ of edges, where each edge is an unordered
pair of vertices. The vertices $v_i$ and $v_j$ are \emph{adjacent}
if $\{v_i,v_j\}$ is an edge. Then $v_j$ is a \emph{neighbour} of
$v_{i}$. All the neighbours of vertex $v_{i}$ in graph $\Gamma$ form
the \emph{neighbourhood} of $v_{i}$, and it is denoted by
$N_{\Gamma}(v_{i})$. The \emph{degree of a vertex} $v$ is the number
of vertices adjacent to $v$. A graph is \emph{regular} of degree $k$
if all vertices have the same degree $k$. For a graph
$\Gamma=(V,E)$, suppose that $V'$ is a nonempty subset of $V$. The
subgraph of $\Gamma$ whose vertex set is $V'$ and whose edge set is
the set of those edges of $\Gamma$ that have both ends in $V'$ is
called \emph{the subgraph of $\Gamma$ induced by $V'$} and is
denoted by $\Gamma[V']$, we say that $\Gamma[V']$ is an induced
subgraph of $\Gamma$. The \emph{adjacency matrix} $A=(a_{ij})$ of a
graph $\Gamma=(V,E)$ is a symmetric $(0,1)$-matrix defined as
follows: $a_{i,j}=1$ if the $i$-th and $j$-th vertices are adjacent,
and $a_{i,j} =0$ otherwise.

\begin{figure}[h,t,b,p]
\begin{center}
\scalebox{0.6}[0.6]{\includegraphics{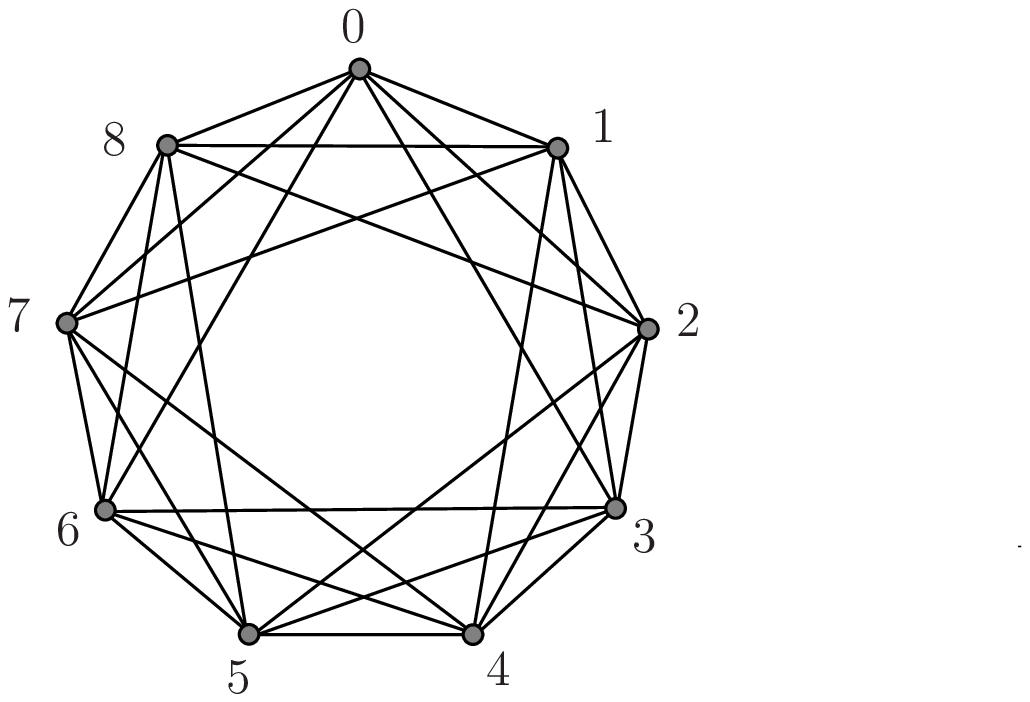}}\\
Figure 1: The circulant graph $C(9,\{1,2,3\})$ \label{fig1}
\end{center}
\end{figure}

Circulant graphs and their various applications are the objects of
intensive study in computer science and discrete mathematics, see
\cite{Bermond, Boesch, Mans, Muzychuk}. Recently, Monakhova
published a survey paper on this subject, see \cite{Monakhova}. Let
$S=\{a_{1},a_{2},\cdots,a_{k}\}$ be a set of integers such that
$0<a_{1}<\cdots<a_{k}<\frac{n+1}{2}$ and let the vertices of an
$n$-vertex graph be labelled $0,1,2,\cdots,n-1$. Then the
\emph{ciculant graph} $C(n,S)$ has $i\pm a_{1},i\pm
a_{2},\cdots,i\pm a_{k} \ (mod\ n)$ adjacent to each vertex $i$. A
\emph{circulant matrix} is obtained by taking an arbitrary first
row, and shifting it cyclically one position to the right in order
to obtain successive rows. We say that a circulant matrix is
generated by its first row. Formally, if the first row of an
$n$-by-$n$ circualant matrix is $a_{0},a_{1},\cdots,a_{n-1}$, then
the $(i,j)^{th}$ element is $a_{j-i}$, where subscripts are taken
modulo $n$. The term circulant graph arises from the fact that the
adjacency matrix for such a graph is a circulant matrix. For
example, Figure \ref{fig1} shows the circulant graph
$C(9,\{1,2,3\})$.

In $2002$, Tonchev \cite{Tonchev} first set up a relationship
between a linear binary code and the adjacency matrix of an
undirected graph. Given a graph on $n$ vertices with adjacency
matrix $A$, one can define two linear codes whose generator matrices
are as follows:

$(a)$ $G=(I;A_n)$,

$(b)$ $G=A_n$,

where $I$ is the identity matrix of order $n$.

The code of type $(a)$ is of length $2n$, dimension $n$, and minimum
Hamming distance $d_n\leq deg_{min}+1$, where $deg_{min}$ is the
minimum degree among the degrees of the vertices in the graph. The
code of type $(b)$ is of length $n$, dimension equal to the rank of
$A$ over the binary field ($2$-rank of $A$), and minimum distance
$d_n\leq deg_{min}$.

Recently, finding optimum codes from graphs has received a wide
attention of many researchers, see \cite{Danielsen, Danielsen2,
Danielsen3, Danielsen4, DP, GGMG, HP, Tonchev, Varbanov}. In
\cite{DP}, Danielsen and Parker showed that two codes are equivalent
if and only if the corresponding graphs are equivalent with respect
to local complementation and graph isomorphism. They used these
facts to classify all codes of lengths up to $12$. In 2012,
Danielsen \cite{Danielsen} focused his attention on additive codes
over $GF(9)$ and transformed the problem of code equivalence into a
problem of graph isomorphism. By an extension technique, they
classify all optimal codes of lengths $11$ and $12$. In fact, a
computer search reveals that circulant graph codes usually contain
many strong codes, and some of these codes have highly regular graph
representations, see \cite{Varbanov}. In \cite{Danielsen}, Danielsen
obtained some optimum additive codes from circulant graphs in 2005.
Later, Varbanov investigated additive circulant graph codes over
$GF(4)$, see \cite{Varbanov}.

In this paper, we introduce a method and find out some optimum
linear codes from circulant graphs. The paper is organized as
follows. In Section $2$, we propose a new method to find linear
optimum codes from circulant graphs. In \cite{Danielsen}, Danielsen
obtained some optimum additive codes from circulant graphs. We get
some optimum linear codes from his result in Section $3$.

\section{New codes from circulant graphs}

We first notice a famous circulant graph, which is called the
$(4,4)$-Ramsey graph. Before introducing this graph, we need some
basic concepts and notions on Ramsey Theory. A \emph{clique} of a
simple graph $\Gamma$ is a subset $S$ of $V$ such that $\Gamma[S]$
is complete. A subset $S$ of $V$ is called \emph{an independent set}
of $\Gamma$ if no two vertices of $S$ are adjacent in $\Gamma$. Let
$r(k,\ell)$ denote the smallest integer such that every graph on
$r(k,\ell)$ vertices contains either a clique of $k$ vertices or an
independent set of $\ell$ vertices. A \emph{$(k,\ell)$-Ramsey graph}
is a graph with $r(k,\ell)-1$ vertices that contains neither a
clique of $k$ vertices nor an independent set of $\ell$ vertices.
The $(4,4)$-Ramsey graph $\Gamma$ (see Figure $2$) is just a
cirulant graph. Let $V(\Gamma)=\{u_1,u_2,\cdots,u_{17}\}$. For the
vertex $u_1$, let
$E_1=\{u_1u_2,u_1u_3,u_1u_5,u_1u_9,u_1u_{10},u_1u_{14},u_1u_{16},u_1u_{17}\}\subseteq
E(\Gamma)$.

\begin{figure}[h,t,b,p]
\begin{center}
\scalebox{0.6}[0.6]{\includegraphics{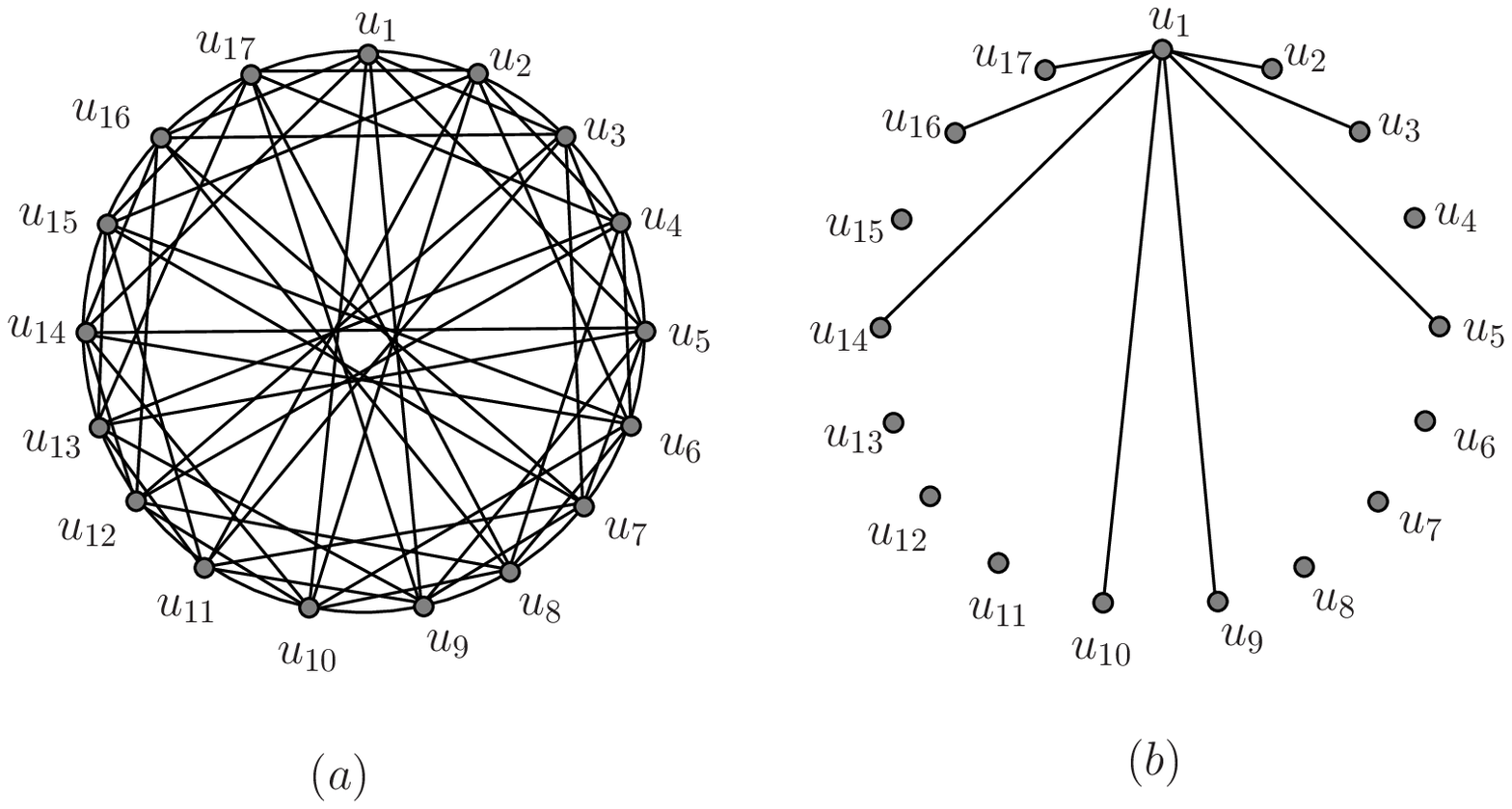}}\\
Figure 2: $(a)$ The $(4,4)$-Ramsey graph $\Gamma$; $(b)$ the edge
induced graph $\Gamma(E_{1})$. \label{fig2}
\end{center}
\end{figure}

The \emph{$(4,4)$-Ramsey graph} is obtained by regarding the
vertices as elements of the field of integers modulo $17$, and
joining two vertices if their difference is a quadratic residue of
$17$ (either $1$, $2$, $4$, $8$, $9$, $13$, $15$ or $16$).  For the
vertex $u_1$, we have
$E_1=\{u_1u_2,u_1u_3,u_1u_5,u_1u_9,u_1u_{10},u_1u_{14},u_1u_{16},$
$u_1u_{17}\}\subseteq E(\Gamma)$. For the vertex $u_2$, we just
rotate the above vertices and edges, that is, we only permit the
existence of the edge set
$E_2=\{u_2u_3,u_2u_4,u_2u_6,u_2u_{10},u_2u_{11},u_2u_{15},$
$u_2u_{17},u_2u_{18}\}\subseteq E(\Gamma)$. For each vertex $u_i\in
V(\Gamma)\setminus \{u_1,u_2\}=\{u_3,u_4,\cdots,u_{17}\}$, we can
also obtain the edge set $E_i \ (3\leq i\leq 17)$. Observe that
$E(\Gamma)=\bigcup_{i=1}^{17}E_i$. For more details, we refer to
\cite{Bondy}. It is clear that the adjacency matrix $A_{17}$ of
$(4,4)$-Ramsey graph is

\begin{equation*}
A_{17}=\left(
\begin{array}{ccccccccccccccccc}
0&1&1&0&1&0&0&0&1&1&0&0&0&1&0&1&1\\
1&0&1&1&0&1&0&0&0&1&1&0&0&0&1&0&1\\
1&1&0&1&1&0&1&0&0&0&1&1&0&0&0&1&0\\
0&1&1&0&1&1&0&1&0&0&0&1&1&0&0&0&1\\
1&0&1&1&0&1&1&0&1&0&0&0&1&1&0&0&0\\
0&1&0&1&1&0&1&1&0&1&0&0&0&1&1&0&0\\
0&0&1&0&1&1&0&1&1&0&1&0&0&0&1&1&0\\
0&0&0&1&0&1&1&0&1&1&0&1&0&0&0&1&1\\
1&0&0&0&1&0&1&1&0&1&1&0&1&0&0&0&1\\
1&1&0&0&0&1&0&1&1&0&1&1&0&1&0&0&0\\
0&1&1&0&0&0&1&0&1&1&0&1&1&0&1&0&0\\
0&0&1&1&0&0&0&1&0&1&1&0&1&1&0&1&0\\
0&0&0&1&1&0&0&0&1&0&1&1&0&1&1&0&1\\
1&0&0&0&1&1&0&0&0&1&0&1&1&0&1&1&0\\
0&1&0&0&0&1&1&0&0&0&1&0&1&1&0&1&1\\
1&0&1&0&0&0&1&1&0&0&0&1&0&1&1&0&1\\
1&1&0&1&0&0&0&1&1&0&0&0&1&0&1&1&0\\
0&1&1&0&1&0&0&0&1&1&0&0&0&1&0&1&1\\
\end{array}\right)
\end{equation*}
It is clear that the adjacency matrix $A_{17}$ of $(4,4)$-Ramsey
graph can be generated by the following vector
$$
\alpha_{17}=(0,1,1,0,1,0,0,0,1,1,0,0,0,1,0,1,1).
$$
Observe that this vector just corresponds to the edge set $E_1$,
which is an expression of the adjacency relation about the vertex
$u_1$. We conclude that the $(4,4)$-Ramsey graph can be determined
by the edge set $E_1$, and the adjacency matrix of this graph is
determined by the above vector. From the relation between a linear
binary code and the adjacency matrix of an undirected graph
introduced by \cite{Tonchev}, we can get an optimum code from the
matrix $(I;A_{17})$.

The above statement suggests the following method.

\textbf{Step 1.} From a circulant graph $\Gamma_n$, we write the
edge set $E_1\subseteq E(\Gamma_n)$, whose elements are incident to
the vertex $u_1\in V(\Gamma_n)$.

\textbf{Step 2.} By the edge set $E_1$, write the vector $\alpha_n$
corresponding to $E_1$.

\textbf{Step 3.} From the vector $\alpha_n$, we generate a circulant
matrix $A_n$.

\textbf{Step 4.} By computer programming, we obtain the minimum
distance of the code $(I;A_{n})$ and determine whether $(I;A_{n})$
is an optimum code.

But, the above method contributes only a few optimum codes. In this
paper, an improved method is introduced by the following statement.

\textbf{Step 1.} From a circulant graph $\Gamma_n$, we write the
edge set $E_1\subseteq E(\Gamma_n)$, whose elements are incident to
the vertex $u_1\in V(\Gamma_n)$.

\textbf{Step 2.} By the edge set $E_1$, write the vector
$$
\alpha_n=(b_1,b_2,\cdots,b_n)
$$
corresponding to $E_1$, where $b_1=0$.

\textbf{Step 3.} From the vector $\alpha_n$, we generate a circulant
matrix $A_n$.

\textbf{Step 4.} Let $L_{n}$ denote the lower bound of the linear
code with length $2n$ from Code Tables. By Algorithm $1$, we obtain
the minimum distance $d_{n}$ of the code $(I;A_{n})$ and determine
whether $d_{n}\geq L_{n}$.

Below is an algorithm (running in SAGE). For more details, we refer
to \cite{Stein}.

\begin{tabbing}
\noindent\rule[0.25\baselineskip]{\textwidth}{2pt}
\\\textbf{Algorithm} 1: Minimum distance of a circulant graph
code\\
\noindent\rule[0.25\baselineskip]{\textwidth}{1pt}\\
Input: the value of $n$, the generator vector $\alpha_n$ of a circulant
graph code $\mathcal{C}_n$\\
Objective: the minimum distance of the circulant graph
code $\mathcal{C}_n$\\
1. input the value of $n$, the generator vector $\alpha_n=(b_1,b_2,\cdots,b_n)$;\\

2. obtain the generator matrix $G=(I;A_{n})$ of the circulant graph
code $\mathcal{C}_n$; \\

3. get the minimum distance of the circulant graph
code $\mathcal{C}_n$.\\[0.2cm]

Take an example, let $n=19$ and
$\alpha_n=(0,1,1,0,1,0,0,0,0,1,1,0,0,0,0,1,0,1,1)$.
 The algorithm\\
details are stated as follows:\\[0.1cm]

\emph{Program}:\\
~~~~~~$n=19$;\\
~~~~~~$a=[0,1,1,0,1,0,0,0,0,1,1,0,0,0,0,1,0,1,1]$\\
~~~~~~$m=matrix(GF(2),[[a[(i-k)\%n]$ for $i$ in $[0..(n-1)]]$ for
$k$ in $[0..n-1]])$;\\
~~~~~~$f=lambda \ s:sum(map(lambda \ x:m[x],s))$;\\
~~~~~~$s=[]$;\\
~~~~~~for $k$ in $[1..8]$:\\
~~~~~~~~~~~$t=min([list(i).count(1)$ for $i$ in Subsets$(range(n),k).map(f)])$;\\
~~~~~~~~~~~$s+=[t]$;\\
~~~~~~~~~~~print $k,t$;\\[0.2cm]

\emph{Output}:~$s_i$ \ \ $s_i'$\\
~~~~~~~~~~~1 \ \ 8\\
~~~~~~~~~~~2 \ \ 6\\
~~~~~~~~~~~3 \ \ 4\\
~~~~~~~~~~~4 \ \ 2\\
~~~~~~~~~~~5 \ \ 2\\
~~~~~~~~~~~6 \ \ 4\\
~~~~~~~~~~~7 \ \ 2\\
~~~~~~~~~~~8 \ \ 2\\

\emph{Result}:~~The elements of the first column
$(s_1,s_2,\cdots,s_8)^{T}$ are the contribution of the matrix $I$ for\\
the weight of a codeword. The elements of the second column
$(s'_1,s'_2,\cdots,s'_8)^{T}$ are the contribution\\ of the matrix
$A_{19}$ for the weight of a codeword. The value of $\min
\{s_i+s_i'\,|\,1\leq i\leq 8\}=6$ is the\\ minimum weight of the
code $\mathcal{C}_{19}$ and then the minimum distance of
the code $\mathcal{C}_{19}$ is also $6$. \\

\noindent\rule[0.25\baselineskip]{\textwidth}{1pt}
\end{tabbing}

\textbf{Step 5.} If the answer is positive, we are done. If not,
i.e., $d_{n}<L_{n}$, then we do the following adjustments of the
elements of the vector $\alpha$. We call codeword $\beta$ a ``bad''
codeword if $wt(\beta)<L_{n}$.

\textbf{Step 5.1}. Find ``bad'' codewords
$\beta_1,\beta_2,\cdots,\beta_{m}$ such that their weights are
$d_n,d_{n}+1,\cdots,d_n+(m-1)$ by Algorithm $2$, where $m=L_n-d_n$.
If there is no codeword with weight $d_n+i \ (0\leq i\leq m-1)$,
then $\beta_{i-1}$ is not under considered.

\textbf{Step 5.2}. For each $\beta_{i} \ (0\leq i\leq m-1)$, we can
find a combination of $\beta_{i}$ by Algorithm $2$. Suppose
$\beta_{i}=\alpha_{n,j_1}+\alpha_{n,j_2}+\cdots +\alpha_{n,j_r}$
where $j_1,j_2,\cdots,j_r\in \{1,2,\cdots,n\}$. Note that
$\alpha_{n,j_k}$ is the $j_k$'s row of the generator matrix. Suppose
$$
\alpha_{n,j_k}=(0,\cdots,0,1,0,\cdots,0\,|\,a_{j_k,1},a_{j_k,2}\cdots,
a_{j_k,n}),
$$
where $a_{j_k,i}=0$ or $a_{j_k,i}=1$ ($1\leq i\leq n$). One can see
that the vector $(a_{j_k,1},a_{j_k,2},\cdots,
a_{j_k,n})=(b_{j_k},b_{j_k+1},\cdots,b_{n},b_1,b_2,\cdots,b_{j_k-1})$.

Below is another algorithm (running in SAGE).

\begin{tabbing}
\noindent\rule[0.25\baselineskip]{\textwidth}{2pt}
\\\textbf{Algorithm} 2: ``bad'' codewords and their combinations\\
\noindent\rule[0.25\baselineskip]{\textwidth}{1pt}\\
Input: the value of $n$, the generator vector $\alpha_n$ of a
circulant
graph code $\mathcal{C}_n$\\
Objective: the minimum distance of the circulant graph
code $\mathcal{C}_n$\\
1. input the value of $n$, the generator vector $\alpha_n=(b_1,b_2,\cdots,b_n)$;\\

2. obtain the generator matrix $G=(I;A_{n})$ of the circulant graph
code $\mathcal{C}_n$; \\

3. get ``bad'' codewords $\beta_1,\beta_2,\cdots,\beta_{m}$ and a
combination of each
$\beta_i \ (1\leq i\leq m)$.\\[0.2cm]

Take an example, let $n=19$ and
$\alpha_n=(0,1,1,0,1,0,0,0,0,1,1,0,0,0,0,1,0,1,1)$.
 The algorithm\\
details are stated as follows:\\[0.1cm]

\emph{Program}:\\
~~~~~~$n=19$;\\
~~~~~~$a=[0,1,1,0,1,0,0,0,0,1,1,0,0,0,0,1,0,1,1]$\\
~~~~~~$m=matrix(GF(2),[[a[(i-k)\%n]$ for $i$ in $[0..(n-1)]]$ for
$k$ in $[0..n-1]])$;\\
~~~~~~$f=lambda \ s:sum(map(lambda \ x:m[x],s))$;\\
~~~~~~$g=lambda \ s:str(sorted(map(lambda \
x:x+1,s))).replace('[','{').replace(']','}')$;\\
~~~~~~$s=[]$;\\
~~~~~~for $k$ in $[1..8]$:\\
~~~~~~~~~~~$t=min([(i,list(f(i)).count(1))$ for $i$ in Subsets$(range(n),k)]$, $key=lambda \ x:x[-1])$;\\
~~~~~~~~~~~$s+=[t]$;\\
~~~~~~~~~~~print $k, t[1], g(t[0])$;\\[0.2cm]

\emph{Output}:~$s_i$ \ \ $s_i'$ \ \ $\{j_1,j_2,\cdots,j_r\}$ \ (as defined in Step 5.2)\\
~~~~~~~~~~~1 \ \ \ 8 \ \ \ \ \{1\}\\
~~~~~~~~~~~2 \ \ \ 6 \ \ \ \ \{1, 9\}\\
~~~~~~~~~~~3 \ \ \ 4 \ \ \ \ \{1, 4, 12\}\\
~~~~~~~~~~~4 \ \ \ 2 \ \ \ \ \{1, 2, 6, 16\}\\
~~~~~~~~~~~5 \ \ \ 2 \ \ \ \ \{1, 2, 8, 11, 14\}\\
~~~~~~~~~~~6 \ \ \ 4 \ \ \ \ \{1, 2, 3, 4, 6, 13\}\\
~~~~~~~~~~~7 \ \ \ 2 \ \ \ \ \{1, 2, 4, 6, 7, 10, 17\}\\
~~~~~~~~~~~8 \ \ \ 2 \ \ \ \ \{1, 2, 3, 6, 7, 8, 12, 16\}\\[0.2cm]
\noindent\rule[0.25\baselineskip]{\textwidth}{1pt}
\end{tabbing}

\textbf{Step 5.3}. Determine whether each element $1$ of the
generator vertex $\alpha_{n}$ is a ``bad'' element in the following
way (Since $b_1=0$, we begin with element $b_2$):

If $b_2=1$, then
$a_{j_1,j_1+1}=a_{j_2,j_2+1}=\cdots=a_{j_r,j_r+1}=b_2=1$. We
calculate the exact value
$$
c_1=\sum_{\ell=1}^{r}a_{j_{\ell},j_{1}+1},\
c_2=\sum_{\ell=1}^{r}a_{j_{\ell},j_{2}+1}, \ \cdots, \
c_{r}=\sum_{\ell=1}^{r}a_{j_{\ell},j_{r}+1}.
$$
Note that $c_i=0$ or $c_i=1$ ($1\leq i\leq r$). Consider the set
$C=\{c_1,c_2,\cdots,c_{r}\}$. If the number of element ``$0$'' in
$C$ is larger than the number of element ``$1$'', then the element
$b_2$ is called \emph{a ``bad'' element} of the generator vector
$\alpha_{n}$. If $b_2$ is a ``bad'' element, then we instead $b_2=1$
by $b'_2=0$ and obtain a new vector
$$
\alpha_n'=(b_1,b'_2,\cdots,b_n)
$$
Then we return to Step $3$. If $b_2$ is not a ``bad'' element or
$b_2=0$, then we consider $b_3$ and continue to determining whether
$b_3$ is a ``bad'' element. The procedure terminates when $b_n$ is
considered.\\

In order to introduce our method clearly, we take the following
example.

Inspired by the above analysis, we hope to consider a circulant
graph of order $19$ having the similar structure with the
$(4,4)$-Ramsey graph.

\textbf{Step 1}. Among all graphs with $19$ vertices, we consider
the graph $\Gamma_{19}$, which can be generated by the edge set
$$
E_1=\{u_1u_2,u_1u_3,u_1u_5,u_1u_{10},u_1u_{11},
u_1u_{16},u_1u_{18},u_1u_{19}\}
$$
Note that this graph has similar structure with the $(4,4)$-Ramsey
graph.

\textbf{Step 2.} By the edge set $E_1$, we write the vector
$$
\alpha_{19}=(0,1,1,0,1,0,0,0,0,1,1,0,0,0,0,1,0,1,1)
$$
corresponding to $E_1$. Obviously,
$b_2=b_{3}=b_{5}=b_{10}=b_{11}=b_{16}=b_{18}=b_{19}=1$. As we see,
$\alpha_{19}$ and $\alpha_{17}$ have a similar distribution of the
elements $0$ and $1$.

\textbf{Step 3}. A circulant matrix can be generated by
$\alpha_{19}$.
\begin{equation*}
A_{19}=\left(
\begin{array}{ccccccccccccccccccc}
0&1&1&0&1&0&0&0&0&1&1&0&0&0&0&1&0&1&1\\
1&0&1&1&0&1&0&0&0&0&1&1&0&0&0&0&1&0&1\\
1&1&0&1&1&0&1&0&0&0&0&1&1&0&0&0&0&1&0\\
0&1&1&0&1&1&0&1&0&0&0&0&1&1&0&0&0&0&1\\
1&0&1&1&0&1&1&0&1&0&0&0&0&1&1&0&0&0&0\\
0&1&0&1&1&0&1&1&0&1&0&0&0&0&1&1&0&0&0\\
0&0&1&0&1&1&0&1&1&0&1&0&0&0&0&1&1&0&0\\
0&0&0&1&0&1&1&0&1&1&0&1&0&0&0&0&1&1&0\\
0&0&0&0&1&0&1&1&0&1&1&0&1&0&0&0&0&1&1\\
1&0&0&0&0&1&0&1&1&0&1&1&0&1&0&0&0&0&1\\
1&1&0&0&0&0&1&0&1&1&0&1&1&0&1&0&0&0&0\\
0&1&1&0&0&0&0&1&0&1&1&0&1&1&0&1&0&0&0\\
0&0&1&1&0&0&0&0&1&0&1&1&0&1&1&0&1&0&0\\
0&0&0&1&1&0&0&0&0&1&0&1&1&0&1&1&0&1&0\\
0&0&0&0&1&1&0&0&0&0&1&0&1&1&0&1&1&0&1\\
1&0&0&0&0&1&1&0&0&0&0&1&0&1&1&0&1&1&0\\
0&1&0&0&0&0&1&1&0&0&0&0&1&0&1&1&0&1&1\\
1&0&1&0&0&0&0&1&1&0&0&0&0&1&0&1&1&0&1\\
1&1&0&1&0&0&0&0&1&1&0&0&0&0&1&0&1&1&0\\
\end{array}\right)
\end{equation*}

\textbf{Step 4.} From the Code Tables, we know that the lower bound
of the minimum distance of linear code $[38,19]$ over $GF(2)$ is
$8$, that is, $L_{19}=8$. By Algorithm $1$, we obtain the minimum
distance $d_{19}$ of the code $(I;A_{19})$ is just $6$, that is,
$d_{19}=6$.

\textbf{Step 5.} Clearly, $6=d_{19}<L_{19}=8$ and
$m=L_{19}-d_{19}=2$.

\textbf{Step 5.1}. From Algorithm $2$, we find two ``bad'' codewords
$$
\beta_1=(1,1,0,0,0,1,0,0,0,0,0,0,0,0,0,1,0,0,0\,|\,0,0,0,0,0,0,0,
1,0,0,0,0,0,1,0,0,0,0,0)
$$
and
$$
\beta_2=(1,0,0,1,0,0,0,0,0,0,0,1,0,0,0,0,0,0,0\,|\,0,1,1,0,
0,1,0,0,0,0,0,0,0,0,0,0,0,1,0)
$$
such that their weights are $6$ and $7$, that is, $wt(\beta_{1})=6$
and $wt(\beta_{2})=7$.

\textbf{Step 5.2}. For $\beta_{1}$, we can find a combination of
$\beta_{1}=\alpha_{19,1}+\alpha_{19,2}+\alpha_{19,6}+\alpha_{19,16}$
by Algorithm $2$, where
$$
\alpha_{19,1}=(1,0,0,0,0,0,0,0,0,0,0,0,0,0,0,0,0,0,0\,|\,0,1,1,0,1,0,
0,0,0,1,1,0,0,0,0,1,0,1,1)
$$
$$
\alpha_{19,2}=(0,1,0,0,0,0,0,0,0,0,0,0,0,0,0,0,0,0,0\,|\,1,0,1,1,0,1,0,0,0,0,1,1,0,0,0,0,1,0,1),
$$
$$
\alpha_{19,6}=(0,0,0,0,0,1,0,0,0,0,0,0,0,0,0,0,0,0,0\,|\,0,1,0,1,1,0,1,1,0,1,0,0,0,0,1,1,0,0,0),
$$
$$
\alpha_{19,16}=(0,0,0,0,0,0,0,0,0,0,0,0,0,0,0,1,0,0,0\,|\,1,0,0,0,0,1,1,0,0,0,0,1,0,1,1,0,1,1,0)
$$
Note that $r=4$, $j_1=1$, $j_2=2$, $j_3=6$ and $j_4=16$.

For $\beta_{2}$, we can find a combination of
$\beta_{2}=\alpha_{19,1}+\alpha_{19,4}+\alpha_{19,12}$ by Algorithm
$2$, where
$$
\alpha_{19,1}=(1,0,0,0,0,0,0,0,0,0,0,0,0,0,0,0,0,0,0\,|\,0,1,1,0,1,0,
0,0,0,1,1,0,0,0,0,1,0,1,1)
$$
$$
\alpha_{19,4}=(0,0,0,1,0,0,0,0,0,0,0,0,0,0,0,0,0,0,0\,|\,0,1,1,0,1,1,0,1,0,0,0,0,1,1,0,0,0,0,1)
$$
$$
\alpha_{19,12}=(0,0,0,0,0,0,0,0,0,0,0,1,0,0,0,0,0,0,0\,|\,0,1,1,0,0,0,0,1,0,1,1,0,1,1,0,1,0,0,0)
$$
Note that $r=3$, $j_1=1$, $j_2=4$ and $j_3=12$.

\textbf{Step 5.3}. Recall that $
\alpha_{19}=(0,1,1,0,1,0,0,0,0,1,1,0,0,0,0,1,0,1,1)$ and
$b_2=b_{3}=b_{5}=b_{10}=b_{11}=b_{16}=b_{18}=b_{19}=1$. Since
$b_2=1$, we consider whether $b_2$ is a ``bad'' element in
$\alpha_{19}$.

For $\beta_{1}$, since $r=4$, $j_1=1$, $j_2=2$, $j_3=6$ and
$j_4=16$, we have
$$
a_{1,2}=a_{2,3}=a_{6,7}=a_{16,17}=b_2=1.
$$
and
\begin{eqnarray*}
c_1&=&\sum_{\ell=1}^{r}a_{j_{\ell},j_{1}+1}=a_{1,2}+a_{2,2}+a_{6,2}+a_{16,2}=0\\
c_2&=&\sum_{\ell=1}^{r}a_{j_{\ell},j_{2}+1}=a_{1,3}+a_{2,3}+a_{6,3}+a_{16,3}=0\\
c_3&=&\sum_{\ell=1}^{r}a_{j_{\ell},j_{3}+1}=a_{1,7}+a_{2,7}+a_{6,7}+a_{16,7}=0\\
c_4&=&\sum_{\ell=1}^{r}a_{j_{\ell},j_{4}+1}=a_{1,17}+a_{2,17}+a_{6,17}+a_{16,17}=0\\
\end{eqnarray*}

For $\beta_{2}$, since $r=3$, $j_1=1$, $j_2=4$ and $j_3=12$, we have
$$
a_{1,2}=a_{4,5}=a_{12,13}=b_2=1.
$$
Then
\begin{eqnarray*}
c_1&=&\sum_{\ell=1}^{r}a_{j_{\ell},j_{1}+1}=a_{1,2}+a_{4,2}+a_{12,2}=1\\
c_2&=&\sum_{\ell=1}^{r}a_{j_{\ell},j_{2}+1}=a_{1,5}+a_{4,5}+a_{12,5}=0\\
c_3&=&\sum_{\ell=1}^{r}a_{j_{\ell},j_{3}+1}=a_{1,13}+a_{4,13}+a_{12,13}=0
\end{eqnarray*}

It is clear that the number of element ``$0$'' in $C$ is larger than
the number of element ``$1$'', then the element $b_2$ is called
\emph{a ``bad'' element} of the generator vector $\alpha_{n}$. We
instead $b_2=1$ by $b_2'=0$ and obtain a new vector
$$
\alpha_{19}'=(0,0,1,0,1,0,0,0,0,1,1,0,0,0,0,1,0,1,1)
$$
Then we return to Step $3$. The circulant matrix $A_{19}'$ generated
by $\alpha_{19}'$ is

\begin{equation*}
A_{19}'=\left(
\begin{array}{ccccccccccccccccccc}
0&0&1&0&1&0&0&0&0&1&1&0&0&0&0&1&0&1&1\\
1&0&0&1&0&1&0&0&0&0&1&1&0&0&0&0&1&0&1\\
1&1&0&0&1&0&1&0&0&0&0&1&1&0&0&0&0&1&0\\
0&1&1&0&0&1&0&1&0&0&0&0&1&1&0&0&0&0&1\\
1&0&1&1&0&0&1&0&1&0&0&0&0&1&1&0&0&0&0\\
0&1&0&1&1&0&0&1&0&1&0&0&0&0&1&1&0&0&0\\
0&0&1&0&1&1&0&0&1&0&1&0&0&0&0&1&1&0&0\\
0&0&0&1&0&1&1&0&0&1&0&1&0&0&0&0&1&1&0\\
0&0&0&0&1&0&1&1&0&0&1&0&1&0&0&0&0&1&1\\
1&0&0&0&0&1&0&1&1&0&0&1&0&1&0&0&0&0&1\\
1&1&0&0&0&0&1&0&1&1&0&0&1&0&1&0&0&0&0\\
0&1&1&0&0&0&0&1&0&1&1&0&0&1&0&1&0&0&0\\
0&0&1&1&0&0&0&0&1&0&1&1&0&0&1&0&1&0&0\\
0&0&0&1&1&0&0&0&0&1&0&1&1&0&0&1&0&1&0\\
0&0&0&0&1&1&0&0&0&0&1&0&1&1&0&0&1&0&1\\
1&0&0&0&0&1&1&0&0&0&0&1&0&1&1&0&0&1&0\\
0&1&0&0&0&0&1&1&0&0&0&0&1&0&1&1&0&0&1\\
1&0&1&0&0&0&0&1&1&0&0&0&0&1&0&1&1&0&0\\
0&1&0&1&0&0&0&0&1&1&0&0&0&0&1&0&1&1&0\\
\end{array}\right)
\end{equation*}

Let us now investigate the linear code $\mathcal {C}_{19}'$ with
generator matrix $G=(I; A_{19}')$. By Algorithm $1$, we get that the
minimum distance $d_{19}'$ of linear code $\mathcal{C}_{19}'$ is
$8$. Thus, the graph code $\mathcal {C}_{19}'$ attains the lower
bound $8$, and hence the code $\mathcal {C}_{19}'$ is a proposed
optimum code over $GF(2)$. The weight enumerator of the code
$\mathcal {C}_{19}'$ is
\begin{eqnarray*}
W_{\mathcal {C}_{19}'}(z)&=&1+133z^{8}+2052z^{10}+10108z^{12}+36575z^{14}+85595z^{16}\\
& &+127680z^{18}+127680z^{20}+85595z^{22}+36575z^{24}+10108z^{26}\\
& &+2052z^{28}+133z^{30}+z^{38}.
\end{eqnarray*}

With the above approach and algorothms, we can also find three other proposed optimum linear codes by the generator matrices $G=(I; A_{19}'')$, $G=(I;A_{19}''')$ and $G=(I; A_{19}'''')$. The circulant matrices $A_{19}''$, $A_{19}'''$ and $A_{19}''''$ are separately generated by
$$
\alpha_{19}''=(0,1,1,0,1,0,0,0,0,0,1,0,0,0,0,1,0,1,1),
$$
$$
\alpha_{19}'''=(0,1,1,0,1,0,0,0,0,1,0,0,0,0,0,1,0,1,1),
$$
$$
\alpha_{19}''''=(0,1,1,0,1,0,0,0,0,1,1,0,0,0,0,1,0,1,0).
$$
The weight enumerator of the code $\mathcal {C}_{19}''$ is
\begin{eqnarray*}
W_{\mathcal {C}_{19}''}(z)&=&1+190z^{8}+1767z^{10}+10507z^{12}+36860z^{14}+84341z^{16}\\
& &+128478z^{18}+128478z^{20}+84341z^{22}+36860z^{24}+10507z^{26}\\
& &+1767z^{28}+190z^{30}+z^{38}
\end{eqnarray*}
One can also check that the weight enumerators of the codes
$\mathcal {C}_{19}'''$ and $\mathcal {C}_{19}''''$ are equal to the
ones of $\mathcal {C}_{19}''$ and $\mathcal {C}_{19}'$,
respectively.

In addition, we consider graphs with large number of vertices with
similar approach. From the Code Tables, we know that the lower and
upper bounds of the minimum distance of binary linear code $[50,25]$
over $GF(2)$ are $10$ and $12$. Let $\Gamma^1_{25}$ and $\Gamma^2_{25}$ be two circulant graphs obtained from the edge sets $$
E^1_1=\{u_1u_2,u_1u_3,u_1u_5,u_1u_7, u_1u_{12},
u_1u_{13},u_1u_{18},u_1u_{20},u_1u_{22},u_1u_{24},u_1u_{25}\}\subseteq E(\Gamma^1_{25})
$$
and
$$
E^2_1=\{v_1v_2,v_1v_3,u_1v_5,v_1v_7,v_1v_{9},
v_1v_{14}, v_1v_{15},v_1v_{20},v_1v_{22},v_1v_{24},v_1v_{25},\}\subseteq E(\Gamma^2_{25}),
$$
respectively, in which $u_1\in V(\Gamma^1_{25})$ and $v_1\in V(\Gamma^2_{25})$. By the edge sets $E^1_1$ and $E^2_1$, we write the generator vectors
$$
\alpha^1_{25}=(0,1,1,0,1,0,1,0,0,0,0,1,1,0,0,0,0,1,0,1,0,1,0,1,1)
$$
and
$$
\alpha^2_{25}=(0,1,1,0,1,0,1,0,1,0,0,0,0,1,1,0,0,0,0,1,0,1,0,1,1).
$$
Then by Algorithm $1$ and Algorithm $2$, we find the linear codes
$\mathcal{C}'_{25}$ and $\mathcal{C}_{25}''$ with generator matrices
$G'=(I; A'_{25})$ and $G''=(I; A_{25}'')$, where the
circulant matrices $A'_{25}$ and $A_{25}''$ are generated
by
$$
\alpha_{25}'=(0,1,0,0,1,0,1,0,0,0,0,1,1,0,0,0,0,1,0,1,0,1,0,1,1)
$$
$$
\alpha_{25}''=(0,1,1,0,1,0,1,0,1,0,0,0,0,1,1,0,0,0,0,1,0,1,0,0,1)
$$
The minimum distance of codes $\mathcal{C}'_{25}$ and
$\mathcal{C}_{25}''$ are $10$. So the graph codes $\mathcal{C}'_{25}$
and $\mathcal{C}_{25}''$ attain the lower bound $10$ and hence the
codes $\mathcal{C}'_{25}$ and $\mathcal{C}_{25}''$ are proposed
optimum codes over $GF(2)$. The weight enumerator of the code
$\mathcal {C}_{25}'$ is
\begin{eqnarray*}
W_{\mathcal
{C}_{25}'}(z)&=&1+225z^{10}+1250z^{11}+3825z^{12}+11525z^{13}+28050z^{14}+64005z^{15}+147075z^{16}\\
& &+294975z^{17}+535075z^{18}+9111100z^{19}+1409205z^{20}+1999925z^{21}+2642200z^{22}\\
& &+3219675z^{23}+3623325z^{24}+377243z^{25}+3621975z^{26}+3216050z^{27}+2643475z^{28}\\
& &+2009175z^{29}+1408010z^{30}+ 904475 z^{31}+
535400z^{32}+292725z^{33}+147525 z^{34}\\
& &+68880z^{35}+27975 z^{36}+9775 z^{37}+ 3500z^{38}+1125 z^{39}+375
z^{40}+125 z^{41}
\end{eqnarray*}
The weight enumerator of the code $\mathcal {C}_{25}''$ is
\begin{eqnarray*}
W_{\mathcal
{C}_{25}''}(z)&=&1+225z^{10}+1250z^{11}+3825z^{12}+11525z^{13}+28050z^{14}+64005z^{15}+147075z^{16}\\
& &+294975z^{17}+535075z^{18}+9111100z^{19}+1409205z^{20}+1999925z^{21}+2642200z^{22}\\
& &+3219675z^{23}+3623325z^{24}+377243z^{25}+3621975z^{26}+3216050z^{27}+2643475z^{28}\\
& &+2009175z^{29}+1408010z^{30}+ 904475 z^{31}+ 535400z^{32}+292725z^{33}+147525 z^{34}\\
& &+68880z^{35}+27975 z^{36}+9775 z^{37}+ 3500z^{38}+1125 z^{39}+375
z^{40}+125 z^{41}
\end{eqnarray*}

\section{New codes from additive codes obtained by Danielsen}

In \cite{Danielsen}, Danielsen got some optimum additive codes. One
of them is the optimum additive code $(30,2^{30},12)$ obtained from
the vector
$$
\beta_{30}=(\omega,0,1,1,0,0,0,0,1,1,0,1,1,1,1,1,1,1,1,1,0,1,1,0,0,0,0,1,1,0),
$$
which corresponds to a circulant graph of order $30$. Change the
``$\omega$'' to ``$0$'', we obtain the following vector.
$$
\alpha_{30}=(0,0,1,1,0,0,0,0,1,1,0,1,1,1,1,1,1,1,1,1,0,1,1,0,0,0,0,1,1,0).
$$
Denote by $A_{30}$ the circulant matrix generated by $\alpha_{30}$.
We now consider the binary linear codes $\mathcal{C}_{30}$ with
generator matrix $G=(I; A_{30})$. By Algorithm $1$, the minimum
distance $d_{30}$ of linear code $\mathcal{C}_{30}$ is $12$. So
$\mathcal{C}_{30}$ is a proposed optimum linear code. The weight
enumerator of the code $\mathcal {C}_{30}$ is
\begin{eqnarray*}
W_{\mathcal {C}_{30}}(z)&=&1+4060z^{12}+24360z^{14}+294930z^{16}+1728400z^{18}+7758400z^{20}\\
& &+26336640z^{22}+67403540z^{24}+129936240z^{26}+192974265 z^{28}+220819632z^{30}\\
& &+192974265z^{32}+129936240z^{34}+67403540z^{36}+26336640
z^{38}+7758660z^{40}\\
& &+1728400z^{42}+294930z^{44}+24360z^{46}+4060
z^{48}+z^{60}\\
\end{eqnarray*}

The above success about obtaining an optimum linear code from
additive codes obtained by Danielsen \cite{Danielsen} gives us more
inspirations since he got more additive codes in \cite{Danielsen}.
In his paper, the code $(15,2^{15})$ from the vector
$$
(\omega,0,1,1,1,0,0,1,1,0,0,1,1,1,0),
$$
which corresponds to a circulant graph of order $15$ is a proposed
optimum additive code. Change the ``$\omega$'' to ``$0$'', we obtain
the following vector.
$$
\alpha_{15}=(0,0,1,1,1,0,0,1,1,0,0,1,1,1,0).
$$
The adjacent matrix generated by $\alpha_{15}$ is

\begin{equation*}
A_{15}=\left(
\begin{array}{ccccccccccccccc}
0&0&1&1&1&0&0&1&1&0&0&1&1&1&0\\
0&0&0&1&1&1&0&0&1&1&0&0&1&1&1\\
1&0&0&0&1&1&1&0&0&1&1&0&0&1&1\\
1&1&0&0&0&1&1&1&0&0&1&1&0&0&1\\
1&1&1&0&0&0&1&1&1&0&0&1&1&0&0\\
0&1&1&1&0&0&0&1&1&1&0&0&1&1&0\\
0&0&1&1&1&0&0&0&1&1&1&0&0&1&1\\
1&0&0&1&1&1&0&0&0&1&1&1&0&0&1\\
1&1&0&0&1&1&1&0&0&0&1&1&1&0&0\\
0&1&1&0&0&1&1&1&0&0&0&1&1&1&0\\
0&0&1&1&0&0&1&1&1&0&0&0&1&1&1\\
1&0&0&1&1&0&0&1&1&1&0&0&0&1&1\\
1&1&0&0&1&1&0&0&1&1&1&0&0&0&1\\
1&1&1&0&0&1&1&0&0&1&1&1&0&0&0\\
0&1&1&1&0&0&1&1&0&0&1&1&1&0&0\\
\end{array}\right)
\end{equation*}
Consider the binary linear code $\mathcal{C}_{15}$ with generator
matrix $G=(I; A_{15})$. Unfortunately, $\mathcal{C}_{15}$ is not a
(proposed) optimum code. So some adjustments are needed for the
elements of $\alpha_{15}$. Applying Algorithm $2$, we obtain a new
vector
$$
\alpha_{15}'=(0,0,1,1,1,0,0,1,1,0,0,1,1,0,0).
$$
The circulant matrix $A_{15}'$ generated by $\alpha_{15}'$ is
\begin{equation*}
A_{15}'=\left(
\begin{array}{ccccccccccccccc}
0&0&1&1&1&0&0&1&1&0&0&1&1&0&0\\
0&0&0&1&1&1&0&0&1&1&0&0&1&1&0\\
0&0&0&0&1&1&1&0&0&1&1&0&0&1&1\\
1&0&0&0&0&1&1&1&0&0&1&1&0&0&1\\
1&1&0&0&0&0&1&1&1&0&0&1&1&0&0\\
0&1&1&0&0&0&0&1&1&1&0&0&1&1&0\\
0&0&1&1&0&0&0&0&1&1&1&0&0&1&1\\
1&0&0&1&1&0&0&0&0&1&1&1&0&0&1\\
1&1&0&0&1&1&0&0&0&0&1&1&1&0&0\\
0&1&1&0&0&1&1&0&0&0&0&1&1&1&0\\
0&0&1&1&0&0&1&1&0&0&0&0&1&1&1\\
1&0&0&1&1&0&0&1&1&0&0&0&0&1&1\\
1&1&0&0&1&1&0&0&1&1&0&0&0&0&1\\
1&1&1&0&0&1&1&0&0&1&1&0&0&0&0\\
0&1&1&1&0&0&1&1&0&0&1&1&0&0&0\\
\end{array}\right)
\end{equation*}
Furthermore, we consider the binary linear code $\mathcal{C}_{15}'$
with generator matrix $G=(I; A_{15}')$.

By Algorithm $1$, the minimum distance $d_{15}'$ of linear code
$\mathcal{C}_{15}'$ is $8$. Again from the Code Tables, we know that
the minimum distance of binary linear code $[30,15]$ over $GF(2)$ is
$8$.  So $\mathcal{C}_{15}'$ is an optimum binary linear code. The
weight enumerator of the code $\mathcal {C}_{15}'$ is
\begin{eqnarray*}
W_{\mathcal
{C}_{15}'}(z)&=&1+450z^{8}+1848z^{10}+5040z^{12}+9045z^{14}+9045z^{16}\\
& &+5040z^{18}+1848z^{20}+450z^{22}+z^{30}.
\end{eqnarray*}


\begin{thebibliography}{11}

\bibitem{Bermond} J. C. Bermond, F. Comellas, D. F. Hsu,
\emph{Distributed loop computer networks: A survey}, J. Parallel
Distributed Comput. 24 (1995), 2-10.

\bibitem{Boesch} F. T. Boesch, J. F. Wang, \emph{Reliable circulant networks with
minimum transmission delay}, IEEE Trans. Circuits Syst. 32 (1985),
1286-1291.



\bibitem{Bondy} J. A. Bondy, U. S. R. Murty,
{\it Graph Theory}, GTM 244, Springer, 2008.


\bibitem{Danielsen} L. E. Danielsen, \emph{On Self-Dual Quantum Codes},
Graphs and Boolean Functions, 2005.

\bibitem{Danielsen2} L. E. Danielsen, \emph{Graph-based classification of self-dual
additive codes over finite field}, Adv. Math. Commun. 3(4) (2009),
329-348.

\bibitem{Danielsen3} L. E. Danielsen, \emph{On the Classification of Hermitian Self-Dual Additive Codes over
GF(9)}, IEEE Trans. Inform. Theory 58(8) (2012), 5500-5511.

\bibitem{Danielsen4} L. E. Danielsen, M. G. Parker, \emph{Directed Graph
Representation of Half-Rate Additive Codes over GF(4)}, Des. Codes
Cryptogr. 59 (2011), 119-130.

\bibitem{DP} L. E. Danielsen, M. G. Parker,
\emph{On the classification of all self-dual additive codes over
$GF(4)$ of length up to $12$}, J. Combin. Theory, Series $A$,
\textbf{113} (2006), 1351-1367.

\bibitem{GGMG} D. G. Glynn, T. A. Gulliver, J. G. Marks, M. K. Gupta,
\emph{The Geometry of Additive Quantum Codes}, Preface, Springer,
2006.

\bibitem{HP} W. C. Huffman, Vera Pless,
{\it Fundamentals of Error-Correcting Codes}, Cambridge University,
2003.

\bibitem{Mans} B. Mans, F. Pappalardi, I. Shparlinski, \emph{On the spectral Adam
property for circulant graphs}, Discrete Math. 254(1-3) (2002),
309-329.

\bibitem{Meijer} P. T. Meijer,
{\it Connectivities and Diameters of Circulant Graps}, B. Sc.
(Honors), Simon Fraser University, 1987.

\bibitem{Monakhova} E. A. Monakhova, \emph{A survey on undirected
circulant graphs}, Discrete Mathematics, Algorithms and
Applications, 4 (1) (2012), DOI: 10.1142/S1793830912500024.

\bibitem{Muzychuk} M. E. Muzychuk, G. Tinhofer, \emph{Recognizing
circulant graphs of prime order in polynomial time}, Electron. J.
Combin. 5(1) (1998), 501-528.

\bibitem{Stein} W. A. Stein et al., \emph{Sage Mathematics Software (Version
  6.1.1)}, The Sage Development Team, 2014, http://www.sagemath.org.

\bibitem{Tonchev} V. Tonchev,
\emph{Error-correcting codes from graphs}, Disrete Math.,
\textbf{257} (2002), 549-557.

\bibitem{Varbanov} Z. Varbanov,
\emph{Additive circulent graph codes over GF(4)}, Math. Maced. 6
(2008), 73-79.

\end{thebibliography}
\end{document}